\documentclass[12pt]{article}
\usepackage{amssymb}

\newcommand{\qed}{\begin{flushright} \vspace{-1pc} $\square$
\end{flushright}}
\newenvironment{dfn}{\bigskip \noindent \bf Definition \rm}{\bigskip}
\newenvironment{rem}{\bigskip \noindent \bf Remark \rm}{\bigskip}
\newenvironment{example}{\bigskip \noindent \bf Example \rm}{\bigskip}
\newenvironment{proof}{\bigskip \noindent \bf Proof: \rm}{\qed \bigskip}
\newtheorem{thm}{Theorem}
\newtheorem{cor}[thm]{Corollary}
\newtheorem{lemma}[thm]{Lemma}
\newfont{\knot}{mini scaled 1000}
\begin{document}
\begin{center}
{\Large\bf On a Planarity Criterion \\[1ex]
Coming from Knot Theory} \\[3ex]
\large J\"org Sawollek\footnote{Fachbereich Mathematik, Universit\"at
Dortmund, 44221 Dortmund, Germany \\ {\em E-mail:\/}
sawollek@math.uni-dortmund.de \\ {\em WWW:\/}
http://www.mathematik.uni-dortmund.de/lsv/sawollek} \\[1ex]
January 25, 2000 (revised: August 29, 2002)
\end{center}
\vspace{4ex}

\begin{abstract}
A graph $G$ is called {\em minimalizable\/} if a diagram with minimal
crossing number can be obtained from an arbitrary diagram of $G$ by
crossing changes. If, furthermore, the minimal diagram is unique up to
crossing changes then $G$ is called {\em strongly minimalizable}. In
this article, it is explained how minimalizability of a graph is related
to its automorphism group and it is shown that a graph is strongly
minimalizable if the automorphism group is trivial or isomorphic to a
product of symmetric groups. Then, the treatment of crossing number
problems in graph theory by knot theoretical means is discussed and, as
an example, a planarity criterion for minimalizable graphs is given.
\\[1ex]
\\[1ex]
\end{abstract}

\section{Introduction}

One of the hardest tasks in topological graph theory is the problem to
determine the minimal crossing number of graph embeddings. Only very few
results are known and there are several outstanding conjectures
concerning the crossing numbers of, e.g., the complete graph $K_n$, the
complete bipartite graph $K_{n,m}$, and the product of cycles $C_n
\times C_m$ for arbitrary values of $m, n \in \mathbb N$.

A main obstacle for finding lower bounds of crossing numbers by means of
knot theory is the fact that, in contrast to knot diagrams, the diagram
of a planar graph in general cannot be unknotted by crossing changes.
Likewise, it is in general not possible to obtain a diagram with minimal
crossing number by crossing changes in an arbitrary diagram. In Section
\ref{mini} of this paper, the class of graphs having this property is
investigated. A graph belonging to this class is called {\em
minimalizable\/} and, in the case that a minimalizable graph possesses a
unique minimal diagram up to crossing changes, it is called {\em
strongly minimalizable}. The r\^ole that the graph's automorphism group
is playing with respect to being minimalizable is investigated and
examples of strongly minimalizable graphs are given, namely, those
graphs whose automorphism group is either trivial or an appropriate
product of symmetric groups.

In Section \ref{crossing}, a general method is described how crossing
number problems in graph theory can be reduced to determine crossing
numbers of finitely many graph diagrams by knot theoretical means. As an
example, a planarity criterion for minimalizable graphs that arises from
knot theory is stated and it is proved in Section \ref{proof}.

\section{Minimalizable Graphs}
\label{mini}

The graphs considered in the following are allowed to have multiple
edges and loops. A {\em topological graph\/} is a 1-dimensional cell
complex which is related to an abstract graph in the obvious way. If $G$
is a topological graph, then a {\em graph $\mathcal G$ in $\mathbb
R^3$\/} is the image of an embedding of $G$ into $\mathbb R^3$. Two
graphs ${\mathcal G}_1$, ${\mathcal G}_2$ in $\mathbb R^3$ are called
{\em equivalent\/} or {\em ambient isotopic\/} if there exists an
orientation preserving autohomeomorphism of $\mathbb R^3$ which maps
${\mathcal G}_1$ onto ${\mathcal G}_2$. Embeddings of topological graphs
in $\mathbb R^3$ can be examined via {\em graph diagrams}, i.e., images
under regular projections to an appropriate plane equipped with
over-under information at double points. Two graph diagrams $D$ and $D'$
are called {\em equivalent\/} or {\em ambient isotopic\/} if the
corresponding graphs in $\mathbb R^3$ are equivalent. Equivalent graph
diagrams can be transformed into each other by a finite sequence of
so-called {\em Reidemeister moves\/} combined with orientation
preserving homeomorphisms of the plane to itself, see \cite{kauf89} or
\cite{yet}.

The {\em crossing number cr(G)\/} of a graph $G$ is the minimal number
of double points in a regular projection of any graph embedding in
3-space. In graph theory, the equivalent definition of $cr(G)$ as
minimal number of crossings in a {\em good drawing\/} of $G$ is more
common, see \cite{erdguy} for an introduction.

\begin{dfn}
A graph $G$ is called {\em minimalizable \/} if a diagram with minimal
crossing number can be obtained from an arbitrary diagram of $G$ by a
choice of crossing changes followed by an ambient isotopy. $G$ is called
{\em strongly minimalizable\/} if it is minimalizable and possesses a
minimal diagram that is unique up to crossing changes followed by an
ambient isotopy.
\end{dfn}

\begin{rem}
It is well-known that the (planar) graph with one vertex and one edge
corresponding to a classical knot is strongly minimalizable. There exist
graphs, even planar ones, that are not minimalizable, see \cite{tan}.
\end{rem}

\begin{thm}
\label{altdef}
Any two diagrams of a strongly minimalizable graph are equivalent up
to crossing changes, i.e., they can be transformed into each other by a
finite sequence of crossing changes and ambient isotopies.
\qed
\end{thm}

\begin{rem}
Observe that every finite sequence of crossing changes and ambient
isotopies applied to a graph diagram can be reduced such that all
crossing changes are realized in a single diagram. This can be seen in
the same way as for the equivalent definitions of the unknotting number
of a knot via crossing changes and ambient isotopies, see \cite{adams}.
\end{rem}

\begin{thm}
\label{triv}
A planar minimalizable graph is strongly minimalizable.
\end{thm}

\begin{proof}
This follows immediately from a result by Lipson \cite{lip}, Corollary
6.
\end{proof}

\begin{rem}
In \cite{tan}, a planar minimalizable graph is called {\em
trivializable\/} and it is shown that a subgraph of a trivializable
graph again is trivializable. In general, a subgraph of a strongly
minimalizable graph is not strongly minimalizable, nor even
minimalizable. This follows from the fact that every graph is isomorphic
to a subgraph of the complete graph $K_n$ for some $n \in \mathbb N$ and
this is strongly minimalizable (see Corollary \ref{examplecor} below).
On the other hand, there are infinitely many examples of graphs that are
not minimalizable, see \cite{tan}, \cite{tatsu}, \cite{tsu}.
\end{rem}

As mentioned above, a knot represented by a knot diagram always can be
unknotted by appropriate crossing changes. Likewise, a knotted arc that
connects two points in the plane can be unknotted by crossing changes.
This can be achieved by traveling along the arc from one end to the
other and choosing self-crossings such that each crossing is passed as
an overcrossing when reached for the first time.

Furthermore, up to crossing changes followed by an ambient isotopy,
there are only finitely many ways of drawing a graph in the plane. The
number of essentially different drawings depends on the graph's
automorphism group.

\begin{lemma}
\label{unique}
Let $G$ be a graph with vertices $v_1, \ldots, v_n$ and let $w_1,
\ldots, w_n$ be distinct points in the plane. Then, up to crossing
changes followed by an ambient isotopy, $G$ has a unique diagram such
that the vertex $v_i$ corresponds to the point $w_i$ for $i = 1, \ldots,
n$.
\end{lemma}

\begin{proof}
The proof is carried out by induction on the number $k \geq 0$ of graph
edges. If $G$ has no edge then there is nothing to show since the $n$
points in the plane can be arranged arbitrarily by ambient isotopy.

Now let $G$ have $k \geq 1$ edges and let $D$ and $D'$ be diagrams of
$G$ such that $w_1, \ldots, w_n$ correspond to the graph vertices.
Choose an arbitrary edge $e$ of $G$. By induction hypothesis, the
diagrams arising from $D$ and $D'$ by deleting the arcs corresponding to
$e$ allow a choice of crossing changes such that the resulting diagrams
are equivalent. Apply these crossing changes to the diagrams $D$ and
$D'$, respectively, and change those crossings in which the arcs $a$ in
$D$ and $a'$ in $D'$ related to $e$ are involved as follows. Choose the
crossings of $a$ with an arc different from $a$ such that $a$ is always
above the other arc and change the self-crossings of $a$ such that $a$
is transformed into an unknotted arc. Apply the same procedure to the
arc $a'$ in $D'$. The two resulting diagrams of $G$ are easily seen to
be equivalent.
\end{proof}

\begin{thm}
\label{identity}
A graph $G$ with trivial automorphism group is strongly minimalizable.
\end{thm}

\begin{proof}
Let $G$ have $n$ vertices. Since $n$ distinct points in the plane can be
moved arbitrarily by ambient isotopy, an embedding of $G$ into the plane
is determined by connecting $n$ fixed points by arcs corresponding to
the incidence relation given by $G$. By Lemma \ref{unique}, the only
ambiguity, up to crossing changes followed by an ambient isotopy, to
connect the points in the plane by arcs arises from the automorphisms of
$G$.
\end{proof}

\begin{rem}
If a graph $G$ possesses no vertices of degree two then there is an
appropriate subdivsion of $G$ that has trivial automorphism group. The
additional vertices are topologically uninteresting but, of course,
important for the property of being minimalizable. Adding vertices of
degree two in the described way has the same effect as colouring the
edges of the graph and considering only graph automorphisms and
modifications of graph diagrams that respect colourings.
\end{rem}

\begin{thm}
\label{group}
Let $G$ be a graph with $n$ vertices. If the automorphism group of $G$
is isomorphic to $S_{n_1} \times \ldots \times S_{n_k}$ with $n_1 +
\ldots + n_k = n$ then $G$ is strongly minimalizable.
\end{thm}

\begin{proof}
As explained in the proof of Theorem \ref{identity}, only the effect of
connecting a fixed set of points by arcs in two different ways arising
from an automorphism of $G$ has to be investigated. Because of the
structure of the automorphism group given, the $n$ points can be
partitioned into subsets corresponding to the partition $n = n_1 +
\ldots + n_k$ and such that whenever two points belonging to different
subsets have to be connected by an arc then every point of the first
subset has to be connected with every point of the second subset and
vice versa. Likewise, whenever two different points belonging to the
same subset have to be connected then any two different points of this
subset have to be connected by an arc and if there is a vertex incident
with a loop then every vertex of the subset is incident with a loop.
Corresponding statements hold in the case of multiple edges. Thus, the
arising diagram of $G$ is unique up to crossing changes followed by an
ambient isotopy.
\end{proof}

\begin{rem}
The construction that is given in the proof of Theorem \ref{group} shows
that the graph $G$ has $K_{n_i,n_j}$ as a subgraph if two points of the
subsets corresponding to $n_i$ and $n_j$ are connected by an arc, and it
has $K_{n_i}$ as a subgraph if two points of the subset corresponding to
$n_i$ are connected by an arc.
\end{rem}

\begin{cor}
\label{examplecor}
For $n, n_1, \ldots, n_k \in \mathbb N$, the complete graph $K_n$ and
the complete $k$-partite graph $K_{n_1, \ldots, n_k}$ are strongly
minimalizable.
\qed
\end{cor}

\section{Crossing Number Problems}
\label{crossing}

Applying the results of the previous section, an algorithm can be given
to reduce the problem of finding the crossing number of a given graph to
the problem of determining the crossing numbers of finitely many graph
diagrams up to ambient isotopy. Start with an arbitrary diagram of the
graph and calculate the crossing numbers of the graph embeddings
corresponding to all possible choices of crossing information for the
diagram's double points. Then change the diagram by connecting vertices
in the plane by arcs in a different manner corresponding to an
automorphism of the graph. As before, calculate the crossing numbers of
all related graphs in $\mathbb R^3$, and carry on until all graph
automorphisms have been applied to the starting diagram. The desired
crossing number of the graph is the minimum taken over all crossing
numbers calculated.

By this procedure, a complicated problem of topological graph theory is
transformed into finitely many knot theoretical problems. Of course, the
determination of crossing numbers for graphs in $\mathbb R^3$ is, in
general, no easy task either. Some results concerning the crossing
numbers of particular classes of graphs can be found in \cite{fourreg}
and \cite{alt}.

The planarity criterion that is stated in the following theorem is an
example how results on crossing numbers for graphs can be deduced from
the calculation of crossing numbers of graph embeddings in $\mathbb
R^3$. A proof of the theorem is given in Section \ref{proof} where the
knot theoretical ingredients are described, too.

\begin{thm}
\label{planar}
Let $G$ be a minimalizable graph and let $D$ be an arbitrary diagram of
$G$. Furthermore, let there exist a graph vertex $v$ in $D$ of degree
four such that the following conditions are fulfilled.
\begin{enumerate}
\item[i)] For each of the four pairs of neighbouring edges of $v$, there
exists a cycle that contains the two edges such that the two cycles in
$D$ belonging to opposite pairs only meet in $v$.
\item[ii)] If $D'$ arises from $D$ by arbitrary crossing changes then
replacing $v$ with an appropriate crossing gives a diagram of an
embedded graph with crossing number at least two.
\end{enumerate}
Then $G$ is non-planar.
\end{thm}

\begin{rem}
Because of the fact that a graph that possesses a non-planar subgraph is
non-planar itself, Theorem \ref{planar} may be applied to graphs that
have no (suitable) vertex of degree four by considering an appropriate
subgraph.

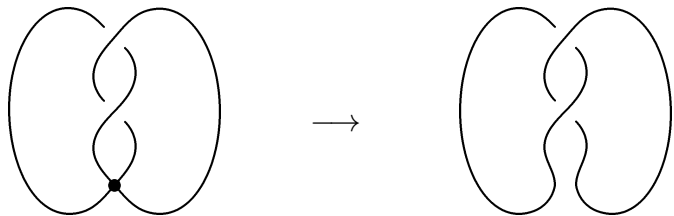
\begin{figure}[htb]
\begin{center}
\begin{picture}(9,2.5)
\put(0,-0.3){\knot o}
\put(4,0.7){$\longrightarrow$}
\put(6,-0.3){\knot p}
\end{picture}
\end{center}
\caption{Forbidden situation yields Hopf link}
\label{linked}
\end{figure}
Observe that the second part of condition i) in Theorem \ref{planar}
excludes the situation depicted in Fig. \ref{linked}, i.e., cutting
through the vertex gives two curves that are linked.
\end{rem}

\begin{example}
\begin{figure}[htb]
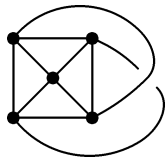

\begin{center}
\knot a
\end{center}
\caption{A minimal diagram of $K_5$}
\label{k5}
\end{figure}
$K_5$ is non-planar as can easily be seen by applying Theorem
\ref{planar} to the diagram depicted in Fig. \ref{k5} for an arbitrary
vertex $v$. Condition i) obviously is fulfilled. For both choices of
over-under information for the diagram's single crossing (indeed, the
two arising diagrams are equivalent), $v$ can be replaced by an
appropriate crossing such that the resulting diagram contains a Hopf
link. Thus condition ii) is fulfilled, too.
\end{example}

Of course, it is well-known that $K_5$ is non-planar and, by
Kuratowski's theorem, that every non-planar graph contains a $K_5$- or a
$K_{3,3}$-minor. The graph considered in the next example contains a
$K_{3,3}$-minor.

\begin{example}
\begin{figure}[htb]
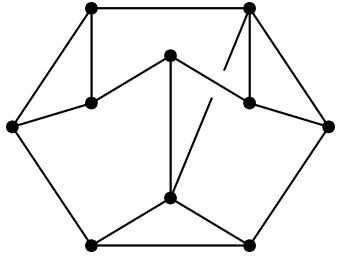

\begin{center}
\knot t
\end{center}
\caption{Non-planar graph}
\label{nonplan}
\end{figure}
With the same argumentation as in the previous example, the graph
depicted in Fig. \ref{nonplan} is non-planar. Again, both of the two
vertices of degree four are appropriate for applying Theorem
\ref{planar}.
\end{example}

In general, very little is known about the behaviour of the crossing
number when a new graph is constructed from one or more given graphs.
For example, it is an open problem wether the crossing number is additive
with respect to a {\em connected sum\/} $G_1 \# G_2$ of two graphs $G_1$
and $G_2$, i.e., two edges $e_1 = \{ v_1, v_2 \} \in G_1$ and $e_2 = \{
w_1, w_2 \} \in G_2$ that are not loops are replaced by edges $e_1' = \{
v_1, w_1 \}$ and $e_2' = \{ v_2, w_2 \}$ (the definition can be extended
to subdivisions of $G_1$ and $G_2$). The corresponding problem for graph
embeddings is well-known in knot theory, and the additivity of crossing
numbers for connected sums of knots and links is an old outstanding
conjecture. In the case of abstract graphs, it is not even clear if a
connected sum of two (strongly) minimalizable graphs is (strongly)
minimalizable.

\begin{thm}
Let $G$ and $G'$ be (strongly) minimalizable graphs. Then the following
hold.
\begin{enumerate}
\item[{\bf a)}] The disjoint union $G \sqcup G'$ of $G$ and $G'$ is
(strongly) minimalizable and
\[ cr(G \sqcup G') \; = \; cr(G) + cr(G') \, . \]
\item[{\bf b)}] Any one-point union $G \bullet G'$ of $G$ and $G'$ is
(strongly) minimalizable and
\[ cr(G \bullet G') \; = \; cr(G) + cr(G') \, . \]
\end{enumerate}
\end{thm}

\begin{proof}
For arbitrary diagrams of $G \sqcup G'$ and $G \bullet G'$,
respectively, a sequence of crossing changes can be chosen such that
arcs corresponding to $G$ always overcross arcs corresponding to $G'$ at
the diagram's double points. Denote the arising diagrams by $D^{\sqcup}$
and $D^{\bullet}$, respectively.

For the graph in $\mathbb R^3$ that belongs to $D^{\sqcup}$, it may be
assumed that the part corresponding to $G$ lies in the half space above
the projection plane and the part corresponding to $G'$ lies beneath the
projection plane. Clearly, there is an equivalent diagram $D \sqcup D'$
which is the disjoint union of a diagram $D$ of $G$ and $D'$ of $G'$.
Since both graphs are minimalizable, there are sequences of crossing
changes that realize $cr(G)$ and $cr(G')$ in $D$ and $D'$, respectively,
showing that $cr(G \sqcup G') \geq cr(G) + cr(G')$. The opposite
inequality holds trivially and it follows that the crossing number is
additive with respect to connected sums. Furthermore, crossing changes
in $D \sqcup D'$ realize $cr(G \sqcup G')$ and therefore $G \sqcup G'$
is (strongly) minimalizable.

Similarly, the graph in $\mathbb R^3$ that belongs to $D^{\bullet}$ may
be thought to consist of a part corresponding to $G$ lying in the upper
half space and a part corresponding to $G'$ lying in the lower half
space except one point $v$ in which the graphs intersect. Deform the
graph corresponding to $G'$ such that its projection is contained
completely in one region of the subdiagram of $D^{\bullet}$
belonging to $G$ (except for the point $v$). This gives an equivalent
diagram $D \bullet D'$ that is a one-point union of a diagram $D$ of $G$
and $D'$ of $G'$. The rest of the proof is completely the same as in the
case of $G \sqcup G'$.
\end{proof}

\section{Proof of Theorem \ref{planar}}
\label{proof}

For the proof of Theorem \ref{planar}, some knot theoretical definitions
and results have to be given. The objects under consideration, as they
are needed here, are described in more detail in \cite{fourreg}. For
general knot theoretical terminology see, e.g., \cite{adams},
\cite{bur}, \cite{kawa}, \cite{lick}, \cite{rol}. In the following, a
{\em link\/} is an embedding of a graph in $\mathbb R^3$ that consists
of one or more disjoint loops.

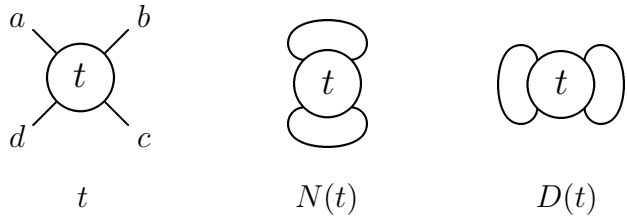
\begin{figure}[hbt]
\begin{center}
\begin{picture}(9,2.5)
\put(0.5,0.7){\knot e}
\put(1.25,1.1){\large $t$}
\put(0.4,1.9){$a$}
\put(2.1,1.9){$b$}
\put(2.1,0.3){$c$}
\put(0.4,0.3){$d$}
\put(1.3,-0.5){$t$}
\put(4.1,0.3){\knot f}
\put(4.55,1){\large $t$}
\put(4.2,-0.5){$N(t)$}
\put(6.9,0.6){\knot g}
\put(7.65,1){\large $t$}
\put(7.4,-0.5){$D(t)$}
\end{picture}
\end{center}
\caption{A tangle and its closures}
\label{tangle}
\end{figure}

A {\em tangle\/} is a part of a link diagram in the form of a disk with
four arcs emerging from it, see Fig. \ref{tangle},
where the tangle's position is indicated by labeling the emerging arcs
with letters $a$, $b$, $c$, $d$ in a clockwise ordering (or simply one
of them with "$a$"). An equivalence relation for tangles is given via
ambient isotopy that fixes the ends of the tangle's arcs. For a tangle
$t$, there are two possible ways to connect the four ends by two arcs
in the plane that do not intersect, see Fig. \ref{tangle}. The arising
link diagrams $N(t)$ and $D(t)$ are called {\em closures\/} of $t$.

\begin{figure}[htb]
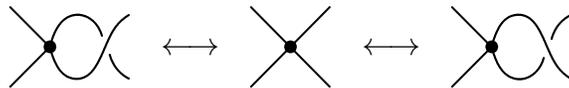

\begin{center}
\begin{tabular}{ccccc}
{\knot b} & $\longleftrightarrow$ & {\knot c} & $\longleftrightarrow$
& {\knot d}
\end{tabular}
\end{center}
\caption{Reidemeister moves of type V}
\label{rfive}
\end{figure}

Denote the two different tangles with crossing number zero by $0$ and
$\infty$, and the two diffent tangles with crossing number one by $1$
and $\overline{1}$. They belong to the class $\cal R$ of {\em rational
tangles\/}. There is an important connection between rational tangles
and {\em Reidemeister moves of type V\/} at a graph vertex as depicted
in Fig. \ref{rfive}.
Indeed, a rational tangle can be defined as the result of applying a
sequence of moves as depicted in Fig. \ref{rfive} to one of the tangles
$0$, $\infty$, $1$, $\overline{1}$ instead of a graph vertex.

\begin{figure}[hbt]
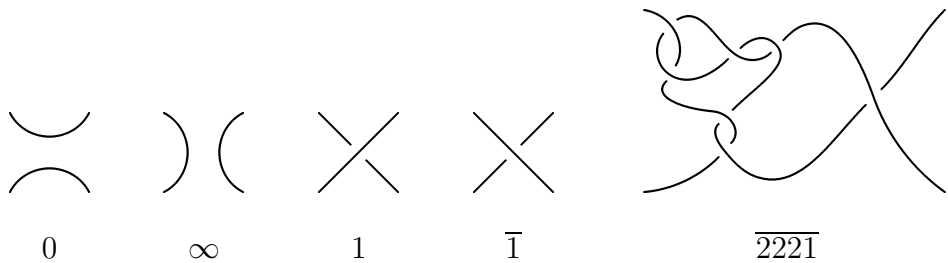

\[ \begin{array}{c@{\hspace{1cm}}c@{\hspace{1cm}}c@{\hspace{1cm}}c@{\hspace{1cm}}c}
\mbox{\knot k} & \mbox{\knot l} & \mbox{\knot m} & \mbox{\knot n} &
\mbox{\knot j} \\[2ex]
0 & \infty & 1 & \overline{1} & \overline{2 2 2 1}
\end{array} \]
\vspace{-4ex}

\caption{Rational tangles in normal form}
\label{rattan}
\end{figure}
Rational tangles can be uniquely classified by rational numbers, see
\cite{con} and \cite{gold}, and there is a normal form for rational
tangles from which the crossing number and an alternating diagram that
realizes this number can easily be read off, see Fig. \ref{rattan}.
For a rational tangle $r$, let $|r|$ denote its crossing number.

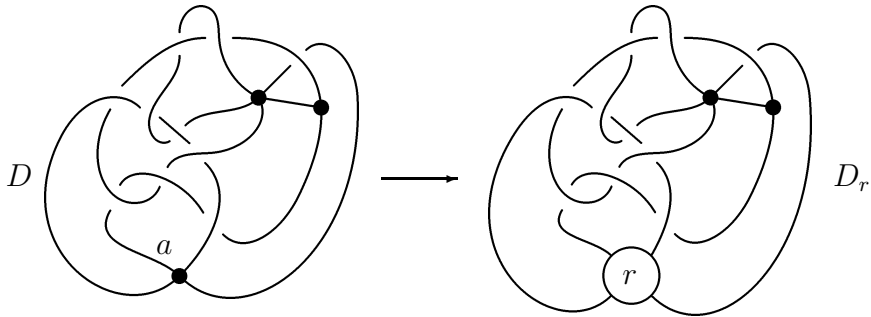
\begin{figure}[htb]
\begin{picture}(12.5,4)
\thinlines
\put(0.5,1.5){$D$}
\put(1.5,0.3){\knot h}
\put(2.5,0.6){$a$}
\put(5.5,1.6){\vector(1,0){1}}
\put(7.5,0.3){\knot i}
\put(8.7,0.2){$r$}
\put(11.5,1.5){$D_r$}
\end{picture}
\caption{Replacing a vertex by a rational tangle}
\label{replace}
\end{figure}

From a given graph diagram $D$ that contains a vertex of degree four,
there can be obtained new graph diagrams by substituting rational
tangles for the graph vertex, see \cite{fourreg}. To do this in a
well-defined way, it is necessary to give an {\em orientation\/} to the
vertex, i.e., labeling an edge incident with the vertex with the letter
$a$. Then a rational tangle can be substituted for the graph vertex in
the obvious way such that the edge labeled with "a" fits together with
the corresponding arc emerging from $r$. An example is given in Fig.
\ref{replace}.

In the same manner as in \cite{fourreg}, where only diagrams of
4-regular graphs are considered, the following theorem can be shown.

\begin{thm}
Let $D$ be a graph diagram that has a vertex of degree four for which an
orientation is chosen. Then, the set $\{ D_r \; | \; r \in {\cal R} \}$
is invariant with respect to ambient isotopy.
\qed
\end{thm}

\begin{figure}[htb]
\begin{center}
\begin{picture}(5,0.5)
\put(1.2,-0.7){\knot q}
\put(1.2,-0.7){\knot r}
\put(3.5,-0.7){\knot r}
\end{picture}
\end{center}
\caption{Non-reduced diagram}
\label{reduced}
\end{figure}
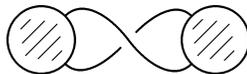

A link diagram $D$ is called {\em reduced\/} if it does not contain a
crossing point $p$ such that $D \setminus \{ p \}$ has more components
than $D$ as depicted in Fig. \ref{reduced}.
Since a reduced alternating link diagram has minimal crossing number,
see \cite{kauf87}, \cite{mura}, \cite{this} for proofs of this famous
{\em Tait Conjecture\/}, the following lemma can readily be deduced from
a rational tangle's normal form.

\begin{lemma}
\label{ratcross}
For every rational tangle $r$ with $|r| \geq 2$, either $N(r)$ or $D(r)$
has crossing number $|r|$.
\qed
\end{lemma}

Using the same technique as for the proof of Lemma \ref{ratcross},
namely, the span of the Jones polynomial which is additive with respect
to connected sums $D_1 \# D_2$ of knot diagrams $D_1$, $D_2$, the
following easily can be shown.

\begin{lemma}
\label{ratsum}
For a rational tangle $r$ with $|r| \geq 2$ and arbitrary knot diagrams
$D_1$, $D_2$, either $D_1 \# N(r) \# D_2$ or $D_1 \# D(r) \# D_2$ has
crossing number $\geq |r|$.
\qed
\end{lemma}

\begin{rem}
Indeed, Lemma \ref{ratsum} holds for any reduced alternating diagram
instead of $N(r)$ or $D(r)$, respectively, and for connected sums with
finitely many diagrams $D_1, \ldots, D_k$. The typical situation that
occurs in the proof of Theorem \ref{planar} is depicted in Fig.
\ref{sum}.
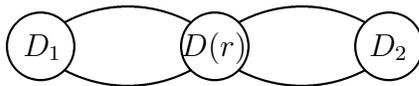
\begin{figure}[htb]
\begin{center}
\begin{picture}(7,0.5)
\put(1.2,-0.7){\knot s}
\put(1.5,-0.3){$D_1$}
\put(3.6,-0.3){$D(r)$}
\put(6.1,-0.3){$D_2$}
\end{picture}
\end{center}
\caption{Connected sum $D_1 \# D(r) \# D_2$}
\label{sum}
\end{figure}
\end{rem}

\noindent
{\bf Proof of Theorem \ref{planar}:}
Assume that $\widetilde{D}$ is a crossing-free diagram of $G$. Following
Theorem \ref{triv}, there exists a diagram $D'$ equivalent with
$\widetilde{D}$ which arises from $D$ by crossing changes. Choose a
vertex-orientation for $v$ in $D'$ and consider the diagrams $D'_r$ with
rational tangles $r$. If $r$ has crossing number at least two then
$D'_r$ likewise has crossing number at least two because of condition
i), which obviously is fulfilled for the diagram $D'$ as well as for
$D$, since a subdiagram of $D'_r$ has crossing number at least two by
Lemma \ref{ratsum}. Furthermore, either $D'_1$ or $D'_{\overline{1}}$
has crossing number at least two because of condition ii). Thus, there
are at most three diagrams $D'_r$, namely, the diagrams $D'_0$,
$D'_{\infty}$, $D'_1$ or the diagrams $D'_0$, $D'_{\infty}$,
$D'_{\overline{1}}$, that have crossing number strictly less than two.

But considering diagrams $\widetilde{D}_r$, corresponding to an arbitrary
vertex of degree four in $\widetilde{D}$ with chosen vertex-orientation,
immediately yields a contradiction since each of the four diagrams
$\widetilde{D}_0$, $\widetilde{D}_{\infty}$, $\widetilde{D}_1$,
$\widetilde{D}_{\overline{1}}$ obviously has crossing number zero or
one.
\qed

\end{document}